\documentstyle[12pt]{article}

\newtheorem{lemma}{Lemma}
[section]

\newtheorem{thm}{Theorem}
[section]
\newtheorem{cor}{Corollary}
[section]
\newtheorem{rmk}{Remark}
[section]

\newcommand{\mysection}[1]{\section{#1}\setcounter{equation}{0}}

\newcommand{\la}{\lambda}

\newcommand{\ptl}{\partial}

\newcommand{\ra}{\rightarrow}
\newcommand{\lra}{\longrightarrow}

\def\a{\alpha}

\def\R2{{\Bbb R}^2}
\def\Bbb{\mathbb}

\def\beq{\begin{equation}}
\def\eeq{\end{equation}}
\def\ba{\begin{array}}
\def\ea{\end{array}}
\def\barr{\begin{eqnarray}}
\def\earr{\end{eqnarray}}

\def\part{\partial}
\def\fr{\frac}

\newcommand{\ls}{\setlength{\baselineskip}{18pt}
                      \setlength{\parskip}{3mm} }

\title {The $C^{\a}$ regularity of a class of hypoelliptic ultraparabolic equations}
\author{ WANG Wendong and ZHANG Liqun
\thanks{ The research is partially supported by the
Chinese NSF under grant 10325104. Email: wangwendong@amss.ac.cn and
lqzhang@math.ac.cn} }

\date{Institute of Mathematics, AMSS, Academia Sinica,
Beijing\\  }

\begin{document}

\maketitle

\begin{abstract}
We obtained the $C^{\a}$ continuity for weak solutions of a class of
ultraparabolic equations with measurable coefficients of the form
${\ptl_t \,u}= \sum_{i,j=1}^{m_0}X_i(a_{ij}(x,t)X_j\,u )+X_0 u$. The
result is proved by simplifying and generalizing our earlier
arguments for the $C^{\a}$ regularity of homogeneous ultraparabolic
equations.
\end{abstract}

{\small keywords: Hypoelliptic, ultraparabolic equations, $C^{\a}$
regularity }

\pagenumbering{arabic}
\mysection{Introduction} \ls \noindent

We are concerned with the regularity of a class of ultraparabolic
equations. One of the typical example of the ultraparabolic
equations is the following equation
$$
\displaystyle \fr{\partial \,u}{\partial \,t}+y \fr{\partial
\,u}{\partial \,x}-u^2\fr{\partial^2 \,u}{\partial
\,y^2}=0,\leqno(1.1)
$$
which is of strong degenerated parabolic type equations, more
precisely, an ultraparabolic type equation. On the other hand, the
equation (1.1), if we consider it as an equation of $\fr 1 u$, has
the divergent form. The recent paper of Pascucci and Polidoro [18],
has proved that the Moser iterative method still works for a class
of ultraparabolic equations with measurable coefficients which are
called homogeneous Kolmogorov-Fokker-Planck equations(or
KFP-equations). By the same technique,  Cinti, Pascucci, Polidoro
[3] consider a nonhomogeneous KFP-equations,  and Cinti, Polidoro
[4] deal with a more general ultraparabolic equation which we will
concentrate on in this paper. Their result shows that for a
non-negative sub-solution $u$ of the ultraparabolic equation they
considered,  the $L^{\infty}$ norm of $u$ is bounded by the $L^p$
norm ($p \ge 1$).

From mathematical points of view, the ultraparabolic equation that
we considered has some special algebraic structures and is
degenerated. There are more and more studies on this problem in
recent years. We have proved that if the weak solution obtained in
[21] of (1.1) is of $C^{\a}$ class, then $u$ is smooth. The second
author [23] has proved $C^{\a}$ property of weak solutions by
Kruzhkov's approach for homogeneous KFP-equations, and the authors
deal with nonhomogeneous KFP-equations in [20]. By simplifying the
cut-off function and generalizing our earlier arguments, we are able
to prove the $C^{\a}$ regularity for weak solutions of more general
ultraparabolic equations. We prove a ${\rm Poincar\acute{e}}$ type
inequality for non-negative weak sub-solutions of (1.2). Then we
apply the inequality to obtain a local priori estimate which implies
the H\"older estimates.

Consider a class of ultraparabolic operator on ${R}^{N+1}$:
$$
\displaystyle  Lu \equiv \sum_{i,j=1}^{m_0}X_i(a_{ij}(x,t)X_j\,u
)+X_0 u - {\ptl_t \,u}=0,\leqno(1.2)
$$
where $(x,t)=z\in { R}^{N+1}$, $1\leq m_0\leq N$, and $X_{j}$'s are
smooth vector fields on $R^N$, for $j= 0,1,\cdots, m_0$.

We follow the notations as in [4]. Let $A=(a_{ij})_{m_0\times m_0},$
and $A_0$ be the identity matrix of $m_0\times m_0.$ Put
$Y=X_0-\ptl_t$,
 and denote
$$ L_1:=\sum_{k=1}^{m_0}X_k^2+Y\leqno(1.3)$$ A curve $\gamma
:[0,T]\rightarrow R^{N+1}$ is called $L_1$-$admissible$, if it is
absolutely continuous and satisfies $$
\gamma'(s)=\sum_{k=1}^{m_0}\la_k
(s)X_k(\gamma(s))+Y(\gamma(s)),\quad \rm{a.e.\,\,in\,[0,T]},$$ for
suitable bounded measurable functions $\la_1(s),\dots,\la_{m_0}(s).$

We make the following assumptions on the operator $L$:

${\bf [H_1]:}$  the coefficients $a_{ij}$, $1\leq i,j\leq m_0$, are
real valued, measurable functions of $(x,t)$. Moreover,
$a_{ij}=a_{ji} \in L^{\infty} ({R}^{N+1})$ and there exists a
$\lambda >0$ such that
$$
\fr{1}{\lambda}\sum_{i=1}^{m_0}\xi_i^2 \leq \sum_{i,j=1}^{m_0}
a_{ij}(x,t)\xi_i \xi_j \leq {\lambda}\sum_{i=1}^{m_0}\xi_i^2
$$
for every $(x,t)\in {R}^{N+1}$, and $\xi \in {R}^{m_0}$.
$X_0=\sum_{i=1}^{N}b_j(x)\ptl_{x_j}$ with smooth functions $b_j(x)$;

${\bf [H_2]:}$ there exists a homogeneous Lie group ${
G}\equiv(R^{N+1},\circ,\delta_{\mu})$ such that \hspace{5pt}(i)
 \hspace{3pt}$X_1,\dots, X_{m_0},Y$ are left translation invariant on
$ {G}$,\\
\quad (ii) \hspace{3pt}$X_1,\dots, X_{m_0}$ are
$\delta_{\mu}$-homogeneous of degree one and Y is
$\delta_{\mu}$-homogeneous of degree two;

 ${\bf [H_3]:}$ for every $(x,t)$, $(\xi, \tau)\in R^{N+1}$
with $t>\tau$, there exists an $L_1$-$admissible$ path $\gamma
:[0,t-\tau]\rightarrow R^{N+1}$ such that $\gamma(0)=(x,t)$,
$\gamma(t-\tau)=(\xi,\tau)$.

 The requirements of ${\bf
[H_2]}$ and ${\bf [H_3]}$ ensure that the operator $L_1$ satisfies
the well-known H\"ormander's hypoellipticity condition by Kogoj and
Lanconelli [8]. We refer to [4] and [8] for more details on the
hypoelliptic type operator on ${R}^{N+1}$.

The Schauder type estimate of (1.2) has been obtained, for example,
Lunardi [12] and Manfredini [14]. Besides, the regularity of weak
solutions have been studied by Bramanti, Cerutti and Manfredini [2],
Polidoro and Ragusa [19], Manfredini and Polidoro [13] assuming a
weak continuity on the coefficient $a_{ij}$. It is quite interesting
whether the weak solution has H\"older regularity under the
assumption ${\bf [H_1]}$. One of the approaches to the H\"older
estimates is to obtain the Harnack type inequality. In the case of
elliptic equations with measurable coefficients, the Harnack
inequality is obtained by J. Moser [15] via an estimate of BMO
functions due to F. John and L. Nirenberg together with the Moser
iteration method. J. Moser [16] also obtained the Harnack inequality
for parabolic equations with measurable coefficients by generalizing
the John-Nirenberg estimates to the parabolic case. Also De Giorgi
developed an approach to obtain the H\"older regularity for elliptic
equations. Another approach to the H\"older estimates is given by S.
N. Kruzhkov [9], [10] based on the Moser iteration to obtain a local
priori estimate, which provides a short proof for the parabolic
equations. Nash [17] introduced another technique relying on the
${\rm Poincar\acute{e}}$ inequality and obtained the H\"older
regularity.

Let $X$ be the gradient with respect to the variables $X_1,
X_2,\cdots, X_{m_0}$, and $Xu=(X_1u, X_2u,\cdots, X_{m_0}u)^T$. We
say that $u$ is a $weak \,solution$ if it satisfies (1.2) in the
distribution sense, that is for any $\phi \in C^{\infty}_0(\Omega)$,
$\Omega$ is a open subset of $R^{N+1}$, then
$$
\int_{\Omega} \phi Yu-(Xu)^T AX\phi = 0, \leqno(1.3)
$$
where $u$, $Xu$, $Yu \in L^2_{\rm loc}(\Omega).$

Our main result is the following theorem.
\begin{thm}
Under the assumptions ${\bf [H_1]}\sim {\bf [H_3]}$, the weak
solution of (1.2) is H\"older continuous.
\end{thm}

\mysection{Some Preliminary and Known Results}

We follow the earlier notations to give some basic properties. For
the more details of the subject, we refer to  Cinti and Polidoro
[4], Kogoj and Lanconelli [8], or Bonfiglioli, Lanconelli and
Uguzzoni[1].

We say a Lie group $G=(R^{N+1},\circ)$ is $homogeneous$ if a family
of dilations $(\delta_{\mu})_{\mu>0}$ exists on G and is an
automorphism of the group: $\delta_{\mu}(z\circ
\zeta)=\delta_{\mu}(z )\circ \delta_{\mu}(\zeta)$, for all
$z,\zeta\in R^{N+1}$ and $\mu>0$, where
$$\delta_{\mu}=diag(\mu^{\a_1},\mu^{\a_2},\dots,\mu^{\a_N},\mu^2),$$
for $1\leq i\leq N$, $\a_i$ is a positive integer, and $\a_1\leq
\a_2\leq\cdots\leq\a_N.$ Moreover, the dilation $\delta_{\mu}$
induces a direct sum decomposition on $R^N$, and $R^N
=V_1\oplus,\cdots,\oplus V_k$.
If we denote $x = x^{(1)} +x^{(2)}
+\cdots+x^{(k)}$ with $ x^{(j )}\in V_j$ , then $$\delta_{\mu}(x, t)
= (D_{\mu} x,\mu^2t),$$ where
$$D_{\mu} (x^{(1)} +x^{(2)} +\cdots+x^{(k)}) = (\mu x^{(1)}
+\mu^2x^{(2)} +\cdots+\mu^k x^{(k)}).$$  Let
$$Q=dimV_1+2dimV_2+\cdots+kdimV_k,$$
then the number $Q+2$ is usually called the $homogeneous$ dimension
of $(R^{N+1},\circ)$ with respect to the dilation $\delta_{\mu}$.

A real function $f(x)$ defined on $R^N$ is called
$\delta_{\mu}$-$homogeneous$ of degree ${m \in R}$, if $f(x)$ does
not vanish identically and, for every $x \in R^N$ and $\mu > 0$, it
holds $$f(\delta_{\mu} (x,0)) = \la^m f(x).$$ A
non-identically-vanishing vector field $X$ is called
$\delta_{\mu}$-$homogeneous$ of degree $m \in R$, if for every
$\phi\in C^{\infty}(R^N)$, $x \in R^N$, and $\mu > 0$, it holds
$$X(\phi(\delta_\mu(x,0))) = \mu^m(X \phi)(\delta_\mu(x,0)).$$

The norm in $R^{N+1}$, related to the group of translations and
dilation to the equation is defined by $$||(x,t)||=r,$$ if $r$ is
the unique positive solution to the equation
$$
\fr{x_1^2}{r^{2\a_1}}+\fr{x_2^2}{r^{2\a_2}}+\cdots+\fr{x_N^2}{r^{2\a_N}}
+\fr{t^2}{r^4}=1,
$$
where $(x,t) \in R^{N+1}\setminus \{0\}$ and by ${\bf [H_2]}$ and
H\"ormander's hypoellipticity condition, we attain
$$
\a_1=\cdots=\a_{m_0}=1, \quad 1<\a_{m_0+1}\leq \cdots\leq \a_N.
$$
And $||(0,0)||=0$. Obviously $$\|\delta_{\mu}(x,t)\|=\mu\|(x,t)\|,$$
for all $(x,t)\in R^{N+1}$. The quasi-distance in $G$ is $$
d(z,\zeta):=\|\zeta^{-1}\circ z\|,\quad \forall z,\zeta\in
R^{N+1},$$ where $$\zeta^{-1}\circ
z=(S(x,t,\xi,\tau),t-\tau)\leqno(2.1)$$ and $S\in R^N$ is smooth
(see [8]). Moreover, for every compact domain $K\in R^{N+1}$, there
exists a positive constant $C_K$ such that $$ C_K^{-1}|z-\zeta|\leq
d(z-\zeta)\leq C_K|z-\zeta|^{\fr 1 k}, \quad \forall z,\zeta\in
K\leqno(2.2)$$ where $|\cdot|$ denotes the usual Euclidean modulus
(see for instance, Prop 11.2 in [7]).

The ball at a point $(x_0,t_0)$ is defined by
$${\cal B}_r(x_0,t_0)=\{(x,t)|\quad ||(x_0,t_0)^{-1}\circ (x,t)||\leq r\},$$
and
$${\cal B}^-_r(x_0,t_0)={\cal B}_r(x_0,t_0)\cap\{t<t_0\}.$$
For convenience, we sometimes use the cube instead of the balls. The
cube at point $(0,0)$ is given by
$$
{\cal C}_r(0,0)=\{(x,t)|\quad |t|\leq r^2,\quad |x_1|\leq r^{\a_1},
\cdots, |x_N|\leq r^{\a_N}\}.
$$
It is easy to see that there exists a constant $\Lambda_1$ such that
$$
{\cal C}_{\fr r {\Lambda_1}}(0,0)\subset{\cal B}_r(0,0)\subset{\cal
C}_{\Lambda_1 r}(0,0),
$$
where $\Lambda_1$ only depends on $Q$ and $N$.

We recall $ L_1=\sum_{k=1}^{m_0}X_k^2+Y, $ whose fundamental
solution $\Gamma_1(\cdot,\zeta)$ with pole in $\zeta\in R^{N+1}$ is
smooth out of the diagonal of $R^{N+1}\times R^{N+1}$, has the
following properties:
$$\ba{lllllll}
(i)\hspace{5pt} \Gamma_1(z,\zeta)=\Gamma_1(\zeta^{-1}\circ
z,0)=\Gamma_1(\zeta^{-1}\circ z), \qquad \forall z, \zeta \in
R^{N+1},\quad z \neq \zeta;
\\ \\
(ii)\hspace{2pt} \Gamma_1(z,\zeta)\geq0,\quad {\rm and}\quad
\Gamma_1(x,t,\xi,\tau)>0\quad{\rm if}\quad t>\tau;\\\\
(iii)\hspace{2pt}\int_{R^N}\Gamma_1(x,t,\xi,\tau)dx=\int_{R^N}\Gamma_1(x,t,\xi,\tau)d\xi=1,\quad{\rm
if}\quad t>\tau;\\\\(iv)\hspace{2pt}\Gamma_1(\delta_{\mu}\circ
z)=\mu^{-Q}\Gamma_1(z),\quad \forall z\neq 0,
\,\mu>0;\ea\leqno(2.3)$$ moreover,
$$
\Gamma_1(z,\zeta)\leq C ||\zeta^{-1}\circ z||^{-Q}, \leqno(2.4)
$$
 for all $z,\zeta\in R^{N+1}$(see [4] or [8]).

A $weak$ $sub$-$solution$ of (1.2) in a domain $\Omega$ is a
function $u$ such that $u$, $Xu$, $Yu \in L^2_{loc}(\Omega)$ and for
any $\phi \in C^{\infty}_0(\Omega)$, $\phi \geq 0$,
$$
\int_{\Omega} \phi Yu-(Xu)^T AX\phi \geq 0. \leqno(2.5)
$$

A result of Cinti and Polidoro obtained by using the Moser's
iterative method (see  Prop 4.4 in [4]) states as follows.

\begin{lemma}
Let $u$ be a non-negative weak sub-solution of (1.2) in $\Omega$.
Let $(x_0,t_0)\in \Omega$ and $\overline{{\cal
B}^-_r(x_0,t_0)}\subset \Omega$ and $p \geq 1$. Then there exists a
positive constant $C$ which depends only on the operator $L$  such
that, for $0 < r\leq 1$
$$
\sup_{{\cal B}^-_{\fr r 2}(x_0,t_0)} u^p \leq \fr
{C}{r^{Q+2}}\int_{{\cal B}^-_r(x_0,t_0)} u^p,\leqno(2.6)
$$
provided that the last integral converges.
\end{lemma}

We make use of a classical potential estimates (see (1.11) in [5])
here to prove the ${\rm Poincar\acute{e}}$ type inequality.

\begin{lemma}
Let $(R^{N+1},\circ)$ is a homogeneous Lie group of homogeneous
dimension $Q+2$, $\a \in (0, Q+2)$ and $G \in C(R^{N+1}\setminus
\{0\})$ be a $\delta_{\mu}$-homogeneous function of degree $\a-Q-2$.
If $f \in L^p(R^{N+1})$ for some $p \in (1,\infty)$, then
$$
G_f(z)\equiv \int_{R^{N+1}} G(\zeta ^{-1}\circ z)f(\zeta)d\zeta,
$$
is defined almost everywhere and there exists a constant
$C=C(Q,p)$ such that
$$
||G_f||_{L^q(R^{N+1})}\leq C \max_{||z||=1} |G(z)|\quad
||f||_{L^p(R^{N+1})},\leqno(2.7)
$$
where $q$ is defined by
$$
\fr 1q =\fr 1p-\fr{\a}{Q+2}.
$$
\end{lemma}
\begin{cor} Let $f\in L^2(R^{N+1})$,  and recall the definitions in [3]
$$
\Gamma_1(f)(z)=\int_{R^{N+1}}\Gamma_1(z,\zeta)f(\zeta) d\zeta,
\qquad \forall z\in R^{N+1},
$$
and
$$
\Gamma_1(X_jf)(z)=-\int_{R^{N+1}}X_j^{(\zeta)}\Gamma_1(z,\zeta)f(\zeta)
d\zeta, \qquad \forall z\in R^{N+1},
$$
where $j=1,\cdots,m_0$, then exists a positive constant $C=C(Q)$
such that
$$
\|\Gamma_1(f)\|_{L^{2\tilde{k}}(R^{N+1})}\leq C\|f\|_{L^2(R^{N+1})},
\leqno(2.8)
$$
and
$$
\|\Gamma_1(X_jf)\|_{L^{2k}(R^{N+1})}\leq C\|f\|_{L^2(R^{N+1})},
\leqno(2.9)
$$
where $\tilde{k}=1+\fr{4}{Q-2}$, $k=1+\fr{2}{Q},$ and
$j=1,\cdots,m_0$.
\end{cor}

\mysection{Proof of Main Theorem}

We may consider the local estimate at a ball centered at $(0,0)$,
since the equation (1.2) is invariant under the left translation
when $a_{ij}$ is constant. The key point in our argument is to
obtain a ${\rm Poincar\acute{e}}$ type inequality. Then by using the
${\rm Poincar\acute{e}}$ type inequality, we prove the following
Lemma 3.5 which is essential in the oscillation estimates in
Kruzhkov's approaches in parabolic case. Then the $C^{\a}$
regularity result follows easily by the standard arguments. We
follow the same route as [23] and [20], but the idea is more simple
and technical. We give them together for completeness.

For convenience, we let $x'=(x_1,\cdots,x_{m_0})$ and $x=(x',
\overline x)$. We consider the estimates in the following cube,
instead of ${\cal B}^-_r$,
$$\ba{lll}
{\cal C}_r^{-}=\{(x,t)| &-r^2\leq t < 0,\, |x'|\leq r,\,
|x_{m_0+1}|\leq (\la' r)^{\a_{m_0+1}}, \cdots, \\ \\ &|x_N|\leq
(\la' r)^{\a_N}\},\ea
$$
where $\la'>1$ is a positive constant, to be decided in (3.8). Let
$$
K_r=\{x'|\quad |x'|\leq r \},
$$
$$
S_r=\{ \overline x\quad|\quad|x_{m_0+1}|\leq (\la' r)^{\a_{m_0+1}},
\cdots, |x_N|\leq (\la' r)^{\a_N}\}.
$$
Let $0<\a, \beta<1$ be constants, for fixed $t$ and $h$, we denote
$$
{\cal N}_{t,h}=\{(x',\overline x)|\quad(x',\overline x)\in K_{\beta
r}\times S_{\beta r},\quad u(\cdot,t) \geq h\}.
$$
By the homogeneousness of $X_j$, $j=1,\cdots,m_0$, we can deduce
$$
X_j=\sum_{i=1}^{m_0}C_i^{(j)}\partial_{x_i}+\sum_{i>m_0}C_i^{(j)}(x)\partial_{x_i},
$$
where $C_i^{(j)}$ is a constant for $i\leq m_0$ and $C_i^{(j)}(x)$
is a polynomial of homogeneous degree $\a_i-1$ for $i> m_0$.
Similarly $$ X_0=\sum_{i>m_0}b_i(x)\partial_{x_i},$$ where $b_i(x)$
is a polynomial of homogeneous degree $\a_i-2$. In the following
discussions, we always assume $r \ll 1$, and  that the constants
$C_i^{(j)}$($i\leq m_0$) and the coefficients of these polynomial
functions $b_j(x)$ and $C_i^{(j)}(x)$($i>m_0$) are bounded by $\la$,
since we can choose $\la$ as a large constant. Moreover, all
constants dependant on $m_0$, $k$ or $Q$ will be denoted by
dependence on $B$.

\begin{lemma}
Suppose that $u(x,t)\geq 0$ be a solution of equation (1.2) in
${\cal C}^-_r$ centered at $(0,0)$ and
$$
mes\{(x,t)\in {\cal C}^-_r, \quad u \geq 1\} \geq \fr 1 2 mes ({\cal
C}^-_r),
$$
then there exist constants $\a$, $\beta$ and $h$, $0<\a, \beta, h<1$
which only depend on $B$, $\la$ and $N$ such that for almost all
$t\in (-\a r^2,0)$,
$$
mes\{{\cal N}_{t,h}\} \geq \fr {1}{11}mes\{ K_{\beta r}\times
S_{\beta r}\}.
$$
\end{lemma}
{\it Proof:} Let
$$
v=\ln^+(\fr{1}{u+h^{\fr 9 8}}),
$$
where $h$ is a constant, $0<h<1$, to be determined later. Then $v$
at points where $v$ is positive, satisfies
$$
\displaystyle \sum_{i,j=1}^{m_0}{X_i(a_{ij}(x,t)X_j\,v )}-(Xv)^TA
Xv+X_0v - {\ptl_t \,v}=0.\leqno(3.1)
$$
Let $\eta(s)$ be a smooth cut-off function so that
$$
\eta(s)=1,\quad \hbox {for} \quad s< \beta r,
$$
$$
\eta(s)=0,\quad \hbox {for} \quad s\geq r.
$$
Moreover, $0\leq\eta \leq 1$ and $|\eta'|\leq \fr {2}{(1-\beta)r}$.

Now we let $\eta_1=\eta(|x'|)$ and $\eta_2=\Pi_{j>m_0}\eta_j$, where
$\eta_j=\eta(\fr1{\la'}|x_j|^{\fr1{\a_j}})$ for $j>m_0$ .

Multiplying $\eta_1^2\eta_2^2$ to (3.1) and integrating by parts on
$K_r\times S_{ r}\times(\tau,t)$
$$
\ba{lllllllll} &&\int_{K_{\beta r}}\int_{S_{ \beta r}}
v(t,x',\overline x)d \overline x dx' +\fr {1}{\la}\int_\tau^t
\int_{K_{ r}}\int_{S_{ r}}\eta_1^2\eta_2^2 \,
|Xv|^2d \overline x dx'dt \\ \\
&\leq& \fr {C(B,\la,N)}
{\beta^{3Q}(1-\beta)^2}(1+\la'^{-1}+\la'^{-2})|S_{\beta
r}|\,|K_{\beta r}|\\\\&&+\int_\tau^t \int_{K_{r}}\int_{S_r}
\eta_1^2\eta_2^2 X_0v d \overline xdx' dt +\int_{K_{r}}\int_{S_{r}}
v(\tau,x',\overline x)d \overline x dx'\\\\&\leq& \fr {C(B,\la,N)}
{\beta^{3Q}(1-\beta)^2}|S_{\beta r}|\,|K_{\beta r}|+\int_\tau^t
\int_{K_{r}}\int_{S_r} \eta_1^2\eta_2^2 X_0v d \overline xdx' dt
\\\\&&+\int_{K_{r}}\int_{S_{r}} v(\tau,x',\overline x)d \overline x
dx',\qquad a.e. \quad\tau, t\in(-r^2,0).\ea\leqno(3.2)
$$
Let $$I_B \equiv \int_{K_{r}}\int_{S_{ r}}\eta_1^2\eta_2^2
\sum_{j>m_0} b_{j}(x)\ptl_{x_j}v d \overline x dx',$$ then
$$
\ba{llllll}  |I_B| &=&|\int_{K_{r}}\int_{S_r}\eta_1^2 \sum_{j>m_0}
(b_{j}(x)\ptl_{x_j}\eta_2^2)v d \overline x dx'|\\\\
&\leq& C(\la, N)\ln(h^{-\fr 9
8})\int_{K_{r}}\int_{S_r}\sum_{j>m_0}|\eta'(|x_j|^{\fr1{\a_j}}
\fr1{\la'})|\,{\fr
1{\la'}}\,|x_j|^{\fr1{\a_j}-1}(\la'r)^{\a_j-2}\\\\
&\leq& {\fr {C(\la ,N)} {(1-\beta)r^2\la'^2}}\beta^{-2Q}|S_{\beta
r}| |K_{\beta r}|\ln(h^{-\fr 9 8}). \ea\leqno (3.4)$$ Integrating by
t to $I_B$, we have
$$
\ba{llllll}\int_\tau^t\int_{K_{r}}\int_{S_{ r}} \eta_1^2\eta_2^2
X_0v d \overline x dx'dt  \leq {\fr {C(\la ,N)}
{(1-\beta)\la'^2}}\beta^{-2Q}|S_{\beta r}| |K_{\beta r}|\ln(h^{-\fr
9 8}). \ea\leqno (3.5)
$$
We shall estimate the measure of the set ${\cal N}_{t,h}$. Let
$$
\mu(t)=mes\{(x',\overline x)|\quad x'\in K_r,\, \overline x \in
S_{r}, \, u(\cdot,t)\geq 1\}.
$$
By our assumption, for $0<\a< \fr 12$
$$
\fr 12 r^2 mes(S_{r})mes(K_{r})\leq \int_{-r^2}^0
\mu(t)dt=\int_{-r^2}^{-\a r^2}\mu(t)dt+\int_{-\a r^2}^{0}\mu(t)dt,
$$
that is
$$
\int_{-r^2}^{-\a r^2}\mu(t)dt\geq (\fr 12-\a)r^2
mes(S_{r})mes(K_{r}),
$$
then there exists a $\tau \in (-r^2,-\a r^2)$, such that
$$
\mu(\tau)\geq (\fr 12-\a)(1-\a)^{-1}
mes(S_{r})mes(K_{r}).\leqno(3.6)
$$
By noticing $v=0$ when $u\geq 1,$ we have
$$
\int_{K_{r}}\int_{S_{ r}} v(\tau,x',\overline x)d \overline x
dx'\leq \fr 12(1-\a)^{-1}mes(S_{r})mes(K_{r})\ln(h^{-\fr 9
8}).\leqno(3.7)
$$
Now we choose $\a$ (near zero), $\beta$ (near one), and $\la'$ large
enough  such that
$$
{\fr {C(\la ,N)} {(1-\beta)\la'^2}}\beta^{-2Q}+\fr{1}{2\beta
^{Q}(1-\a)}\leq \fr 4 5,\leqno(3.8)
$$
and fix them from now on.

By (3.2), (3.5), (3.7) and (3.8), we deduce
$$\ba{lll}
\int_{K_{\beta r}}\int_{S_{\beta r}} v(t,x',\overline x)d \overline
x dx'\\ \\ \leq [\fr {C(B,\la,N)} {\beta^{3Q}(1-\beta)^2} +\fr
45\ln(h^{-\fr 9 8})]mes(K_{\beta r}\times S_{\beta
r}).\ea\leqno(3.9)
$$
When $(x', \bar{x})\notin {\cal N}_{t,h},$, we have
$$\ln(\fr 1 {2h})\leq \ln^+(\fr{1}{h+h^{\fr 9 8}})\leq v,$$
then
$$\ln(\fr 1
{2h})mes(K_{\beta r}\times S_{\beta r}\setminus {\cal N}_{t,h})\leq
\int_{K_{\beta r}}\int_{S_{\beta r}} v(t,x',\overline x)d \overline
x dx'.$$
Since
$$
\fr{C+{\fr 45}\ln(h^{-\fr 98})}{\ln(h^{-1})}\lra \fr
9{10},\qquad\hbox{as} \quad h\ra 0,
$$
then there exists constant $h_1$ such that for $0<h<h_1$ and $t
\in(-\a r^2,0)$
$$
mes(K_{\beta r}\times S_{\beta r}\setminus {\cal N}_{t,h})\leq \fr
{10}{11}mes(K_{\beta r}\times S_{\beta r}).
$$
Then we proved our lemma.

Let $\chi(s)$ be a smooth function given by
$$\ba{ll}
\chi(s)=1 \qquad if \quad s\leq {\theta^{\fr 1 {Q}}} r,\\
\chi(s)=0 \qquad if \quad s> r, \ea
$$
where ${\theta^{\fr 1 {Q}}}<\fr {1}{2}$ is a constant. Moreover, we
assume that
$$
0\leq -\chi'(s) \leq \fr{2}{(1 -{\theta^{\fr 1 {Q}}})r},\quad
|\chi''(s)|\leq \fr{C}{r^2},
$$
and for any $\beta_1, \beta_2,$ with $\theta^{\fr 1
{Q}}<\beta_1<\beta_2<1,$ we have $$|\chi'(s)|\geq
C(\beta_1,\beta_2)r^{-1}>0,$$ if $\beta_1r\leq s\leq \beta_2r.$

For $x\in R^N,$ $t\leq 0$, we set $${\mathcal Q} =\{(x',\bar{x},t)|
-r^2\leq t\leq 0,\, x'\in K_{\fr r \theta}, \,|x_j|\leq({\fr r
{\theta}})^{\a_j}, j=m_0+1,\cdots, N\}.$$ We define the cut off
functions by
$$
\phi_0(x,t)=\chi([\sum_{j=m_0+1}^{N} \fr
{\theta^{2\a_j}x_j^2}{r^{2\a_j-Q}}-C_1 t r^{Q-2} ]^{\fr {1} {Q}}),
$$
$$
\phi_1(x,t)=\chi(\theta |x'|),
$$
$$
\phi(t,x)=\phi_0(t,x)\phi_1(x,t),\leqno(3.10)
$$
where $C_1>1$ is chosen so that
$$
\ba{lllll}  C_1 r^{Q-2}&\geq |\sum_{j>m_0}2
\theta^{2\a_j}b_{j}(x)x_j r^{Q-2\alpha_j}|, \ea
$$
for all $z\in {\mathcal Q}$.

\begin{rmk} By the definition of $\phi$ and the above arguments,
it is easy to check that, for $\theta$,  $r$ small enough and $t\leq 0$\\
(1) $\phi(z)\equiv 1,$ in ${\cal B}^-_{\theta r}$,\\
(2) $\rm{supp}\phi\bigcap \{(x,t)|x\in R^N,t\leq 0\}\subset {\mathcal Q}$,\\
(3) there exists $\a_1>\theta,$ which depends on $C_1,$ such that
$$\{(x,t)| -\a_1r^2\leq t < 0, x'\in K_{\beta r}, \bar{x}\in S_{\beta
r}\}\subseteq \rm{supp}\phi, $$ (4) $0<\phi_0(z)<1,$ for $z\in
\{(x,t)| -\a_1r^2\leq t \leq -\theta r^2, x'\in K_{\beta r},
\bar{x}\in S_{\beta r}\}$.
\end{rmk}

\begin{lemma} Under the above notations, we have
 $$Y \phi_0(z)\leq 0, \quad \rm{for}\quad z\in {\mathcal Q}.$$
\end{lemma}
{\it Proof:} Let $ [\sum_{j=m_0+1}^{N} \fr
{\theta^{2\a_j}x_j^2}{r^{2\a_j-Q}}-C_1 t r^{Q-2}] $ be denoted by
$[\cdots]$. Then
$$\ba{llllllllll} Y \phi_0 &=
\chi'([\cdots]^{\fr{1}{Q}})\fr{1}{Q}[\cdots]^{\fr{1}{Q}-1} [C_1
r^{Q-2} +\sum_{j>m_0}(2 \theta^{2\a_j}b_{j}(x)x_j r^{Q-2\alpha_j})]
 \ea
$$
For the term $b_{j}(x)x_j r^{Q-2\alpha_j}$, since $|b_j(x)|\leq
C(\la,N) ({\fr r \theta})^{\a_j-2}$, we obtain
$$ |\sum_{j>m_0}2 \theta^{2\a_j}b_{j}(x)x_j
r^{Q-2\alpha_j}| \leq C(\la,N) \theta^2 r^{Q-2}.
$$
We can choose a positive constant $C_1>1,$ such that $
C(\la,N)\theta^2<C_1$, then $Y \phi_0(z)\leq 0$ ($z\in {\mathcal
Q}$) holds.

We sometimes abuse the notations of ${\cal B}^-_r$ and ${\cal
C}_r^-$, since there are equivalent. Now we have the following ${\rm
Poincar\acute{e}}$'s type inequality.
\begin{lemma}
Let $w$ be a non-negative weak sub-solution of (1.2) in ${\cal
B}_1^-$. Then there exists a constant $C$, only depends on $B,$
$\la$ and $N$, such that for $r<\theta<1$
$$
\int_{{\cal B}^-_{\theta r}}(w(z)-I_0)_+^2\leq C\theta^2
r^2\int_{{\cal B}^-_{\fr r {\theta}}}|Xw|^2, \leqno(3.11)
$$
where $I_0$ is given by
$$
I_0=sup_{{\cal B}^-_{\theta r}}[I_1(z)+C_2(z)],\leqno(3.12)
$$
and
$$
I_1(z)=\int_{{\cal B}^-_{\fr r {\theta}}}
[-\Gamma_1(z,\cdot)wY\phi](\zeta)d\zeta,\leqno(3.13)
$$
$$
C_2(z)=\int_{{\cal B}^-_{\fr r {\theta}}}
[\sum_{j=1}^{m_0}|X_j^2\phi|\Gamma_1(z,\cdot) w](\zeta)d\zeta,
$$
where $\Gamma_1$ is the fundamental solution,  and $\phi$ is given
by (3.10).
\end{lemma}
 {\it Proof:} We represent $w$ in terms of the fundamental
solution of $\Gamma_1$, i.e.
$$\varphi(z)=-\int_{R_{N+1}}\Gamma_1(z,\zeta)L_1\varphi(\zeta)d\zeta,
 \quad \forall\varphi\in C_0^{\infty}(R^{N+1}).$$ By an approximation and the
support of $\phi$ and $\Gamma_1$,  for $z \in {\cal B}^-_{\theta
r}$, we have
$$\ba{llll}
w(z)&=\int_{{\cal B}^-_{\fr r {\theta}}}  [\langle
A_0X(w\phi),X\Gamma_1(z,\cdot)\rangle
-\Gamma_1(z,\cdot)Y(w\phi)](\zeta)d\zeta \\ \\&=
I_1(z)+I_2(z)+I_3(z)+C_2(z),\ea\leqno(3.14)
$$
where $I_1(z)$ are given by (3.13) and
$$
I_2(z)=\int_{{\cal B}^-_{\fr r {\theta}}} [\langle
(A_0-A)Xw,X\Gamma_1(z,\cdot)\rangle\phi-\Gamma_1(z,\cdot)\langle
(A+A_0)Xw,X\phi\rangle](\zeta)d\zeta,
$$
$$
I_3(z)=\int_{{\cal B}^-_{\fr r {\theta}}} [\langle
AXw,X(\Gamma_1(z,\cdot)\phi)\rangle-\Gamma_1(z,\cdot)\phi
Yw](\zeta)d\zeta.
$$
$$
C_2(z)=\int_{{\cal B}^-_{\fr r {\theta}}} [\langle
A_0X{\phi},X\Gamma_1(z,\cdot)\rangle w+\Gamma_1(z,\cdot)\langle
A_0Xw,X\phi\rangle](\zeta)d\zeta$$
 Note that ${\rm
supp}\phi\bigcap\{\tau\leq 0\}\subset {\mathcal Q}\subset
\overline{{\cal B}^-_{\fr r {\theta}}}$, $z \in {\cal B}^-_{\theta
r}$ and $\langle A_0X{\phi},X\Gamma_1(z,\cdot)\rangle$ vanishes in a
small neighborhood of $z$. Integrating by parts we have
$$
\ba{llll} C_2(z) =\int_{{\cal B}^-_{\fr r {\theta}}}
[\sum_{j=1}^{m_0}|X_j^2\phi|\Gamma_1(z,\cdot) w](\zeta)d\zeta. \ea$$
From our assumption, $w$ is a weak sub-solution of (1.2), and $\phi$
is a test function of this semi-cylinder. In fact, we let
$$
\tilde{\chi}(\tau)=\left\{
\begin{array}{lll} 1\quad &\tau\leq
0,\\1-n\tau\quad &0\leq \tau \leq 1/n,\\ 0\quad &\tau\geq
1/n.\end{array}\right.$$ Then
$\tilde{\chi}(\tau)\phi\Gamma_1(z,\cdot)$ can be a test function
(see [4]). As $n\rightarrow \infty$, we obtain
$\phi\Gamma_1(z,\cdot)$ as a legitimate test function, and
$I_3(z)\leq 0$. Then in ${\cal B}^-_{\theta r}$,
$$
0\leq (w(z)-I_0)_+\leq I_2(z)=I_{21}+I_{22}.
$$
By Corollary 2.1 we have
$$
||I_{21}||_{L^2({\cal B}^-_{\theta r})}\leq C(\la,N)\theta
r||I_{21}||_{L^{2+\fr 4 Q}({\cal B}^-_{\theta r})}\leq C(B, \la, N)
\theta r||Xw||_{L^2({\cal B}^-_{{\fr r {\theta}}})}.\leqno(3.15)
$$
Similarly for $I_{22},$
$$
||I_{22}||_{L^2({\cal B}^-_{\theta r})}\leq |{\cal B}^-_{\theta
r}|^{\fr 12-\fr{Q-2}{2Q+4}} ||I_{22}||_{L^{2\tilde{k}}({\cal
B}^-_{\theta r})}\leq C(B,\la,N) \theta^2 r^2||\,|Xw| \,|X\phi|\,
||_{L^2({\cal B}^-_{{\fr r {\theta}}})},
$$
where $$| X\phi_1|=|\chi'(\theta|\xi'|)\theta X(|\xi'|)|\leq
C(B,\la,N) \fr{\theta}{r}, $$ and $$ | X\phi_0|\leq |\chi'|{\fr 1
Q}[\cdots]^{{\fr 1 Q}-1}\sum_{1\leq i\leq m_0,j>m_0} |\fr
{2C_j^{(i)}(x)\theta^{2\a_j}x_j}{r^{2\a_j-Q}}|\leq
C(B,\la,N)\theta^{\fr 1 Q}r^{-1},
$$ thus
$$||I_{22}||_{L^2({\cal B}^-_{\theta r})}\leq C(B,\la,N) \theta^2
r||Xw||_{L^2({\cal B}^-_{{\fr r {\theta}}})}.$$ Then we proved our
lemma.

Now we apply Lemma 3.3 to the function
$
w= \ln^+\fr{h}{u+h^{\fr 98}}.
$
If $u$ is a weak solution of (1.2), obviously $w$ is a weak
sub-solution. We estimate the value of $I_0$ given by Lemma 3.3.
\begin{lemma}
Under the assumptions of Lemma 3.3, there exist constants
$\lambda_0$, $r_0$ and $r_0<\theta$. $\la_0$ only depends on
constants $\a$, $\beta$, $\lambda$, $B$, $N$, and $\phi$,
$0<\lambda_0<1$, such that for $r<r_0$
$$
|I_0|\leq \lambda_0 \ln(h^{-\fr 1 8}).\leqno(3.16)
$$
\end{lemma}
{\it Proof:} We first come to estimate $C_2(z)$ and as before,
denote $x=(x',\bar{x},t)$ and $\zeta=(\xi',\bar{\xi},\tau)$. Note
 $z \in
{\cal B}^-_{\theta r}$, we have
$$
\ba{llllllllllll} & & |C_2(z)|\\\\
&\leq& \int_{{\cal B}^-_{\fr r {\theta}}}
[\sum_{j=1}^{m_0}|X_j^2\phi|\Gamma_1(z,\cdot) w](\zeta)d\zeta\\\\
&\leq& r^2sup_{\xi\in {\rm
supp}(X\phi)}\sum_{j=1}^{m_0}|X_j^2\phi|\ln (h^{-\fr 1 8}).\quad
(By\,\, $(iii)$ \,\,{\rm in} \,\,(2.3))\ea$$ We only need to
estimate $|X_j^2\phi|$,
$$|X_j^2\phi|\leq |X_j^2\phi_1|+2|X_j\phi_1 X_j\phi_0|+ |X_j^2\phi_0|,$$
where $|X_j\phi_1|=|\theta\chi'(\theta
|\xi'|)\partial_{\xi_j}|\xi'||\leq 2\theta r^{-1},$
$|X_j^2\phi_1|\leq C\theta^{2-\fr1Q}r^{-2} $ and $$
|X_j\phi_0|=|\chi' {\fr 1 Q}[\cdots]^{ {\fr 1
Q}-1}(\sum_{i>m_0}\fr{C^{(j)}_i(\xi)2\xi_i\theta^{2\a_i}}{r^{2\a_i-Q}})|\leq
C(B,\la,N)\theta^{\fr 1 Q}r^{-1},$$ moreover, $$ |X_j^2\phi_0|\leq
C(B,\la,N)\theta^{\fr 1 Q}r^{-2}.
$$Hence
$$|C_2(z)|\leq C(B,\la,N)\theta^{\fr 1 Q}\ln (h^{-\fr 1 8})=C(B,\la,N)\theta^{\a_0} \ln (h^{-\fr 1 8})\leqno(3.17)
$$
where $\a_0={\fr 1 {Q}} >0$.

Since $X_0=\sum_{j>m_0}b_j(x)\partial_{x_j}$, we know $Y\phi=\phi_1
Y\phi_0.$ Now we let $w\equiv 1$, then for $z \in {\cal B}^-_{\theta
r}$ (3.14) gives,
$$
\ba{llll} 1=\int_{{\cal B}^-_{\fr r {\theta}}}[ -
\phi_1\Gamma_1(z,\cdot)Y\phi_0](\zeta)d\zeta +C_2(z)|_{w=1}.\ea\
\leqno(3.18)
$$
 By Lemma 3.2,
$$
-\phi_1\Gamma_1(z,\cdot)Y\phi_0 \geq 0,\leqno(3.19)
$$
we only need to prove $-\phi_1\Gamma_1(z,\cdot)Y\phi_0$ has a
positive lower bound in a domain which $w$ vanishes, and this bound
independent of $\,r$ and small $\theta$. So we can find a $\la_0,$
$0<\la_0<1$, such that this lemma holds and $\la_0$ is independent
of $r$ and small $\theta.$ We observe that the support of $\chi'(s)$
is in the region ${\theta^{\fr 1 {Q}}}r<s< r$, thus for some
${\beta}'< 1$, the set ${\cal B}^-_{{\beta}' r} \setminus {\cal
B}^-_{{\sqrt\theta} r}$ with $\theta r^2/{C_1}\leq |t|\leq \a_1r^2$
is contained in the support of $\phi_1\phi'_0$ . Then we can prove
that the integral of (3.19) on a subset of the domain ${\cal
B}^-_{{\beta}' r} \setminus {\cal B}^-_{{\sqrt\theta} r}$ is lower
bounded by a positive constant.

For $z\in B^-_{\theta r}$, $0<\a_1\leq \a$ and set
$$\zeta\in Z=\{(\xi,\tau)| -\a_1 r^2\leq \tau\leq -\fr {\a_1}
{2}r^2,\,\xi'\in K_{\beta r}, \,\bar{\xi}\in S_{\beta
r},\,w(\xi,\tau)=0\},\leqno(3.20)$$ then $|Z|=C(\a_1,
\beta,B,\la,N)r^{Q+2}$ by Lemma 3.1. We note that when
$\zeta=(\xi,\tau)\in Z $ and $\theta$ is small, $w(\zeta)=0,$
$\phi_1(\zeta)=1$,
$$|\chi'([\cdots]^{\fr{1}{Q}})|\geq C(\a_1,B,\la,N)r^{-1}>0.$$
Consequently
$$
\ba{llllllllllll} \int_Z [-\phi_1\Gamma_1(z,\cdot)Y\phi_0](\zeta)\,d \zeta \\
\\= -\int_Z \phi_1\Gamma_1(z,\cdot)
\chi'([\cdots]^{\fr{1}{Q}})\fr{1}{Q}[\cdots]^{\fr{1}{Q}-1} [C_1
r^{Q-2}
+\sum_{j>m_0}(2\theta^{2\a_j}b_{j}(\xi)\xi_j r^{Q-2\alpha_j})]d\zeta\\ \\
\geq C(B,\la,\a_1,N)\int_Z r^{Q-2}[r^{Q}]^{\fr{1}{Q}-1}
r^{-1}\Gamma_1(\zeta^{-1}\circ z,0)d\zeta \\
\\ \geq C(B,\la,\a_1,N)   \int_Z
r^{-2}\Gamma_1(\zeta^{-1}\circ z,0)d\zeta\\ \\

=C(B,\la,\a,\beta,N)=C_4>0, \ea
$$
where we have used $\Gamma_1(z,\zeta)\geq Cr^{-Q}$, as $\tau \leq
-\fr{\a_1}{2}r^2$ and $z\in B^-_{\theta r}$. In fact, by (iv) in
(2.3) one get
$$\Gamma_1(z,\zeta)=r^{-Q}\Gamma_1(S({x},{t},{\xi},{\tau}),\fr{t-\tau}{r^2}),$$
where $\fr{\a_1}{2}\leq \fr{t-\tau}{r^2}\leq 1$ and S is bounded by
(2.1), hence by the property (ii) in (2.3) of $\Gamma_1$, we have
$\Gamma_1(z,\zeta)\geq C(\a_1)r^{-Q}$. Then we can choose a small
$\theta$ which is fixed from now on and $r_0<\theta$, such that
$$
|I_0|\leq (1-C_4+C_3\theta^{\a_0})\ln(h^{-{\fr 1 8
}})+C_3\theta^{\a_0}\ln(h^{-{\fr 1 8 }})\leq \la_0 \ln(h^{-{\fr 1 8
}})\leqno(3.21)
$$
where $0<r<r_0$, $0<\la_0<1$, depends on $\a$, $\beta$, $B$, $\la$,
$N$, and $\phi$.

\begin{lemma}
Suppose that $u(x,t)\geq 0$ is a solution of equation (1.2) in
${\cal B}^-_r$ centered at $(0,0)$ and $ mes\{(x,t)\in {\cal B}^-_r,
\quad u \geq 1\} \geq \fr 1 2 mes ({\cal B}^-_r). $ Then there exist
constant $\theta$ and $h_0$, $0<\theta, h_0<1$ which only depend on
$B$, $\la$, $\la_0$ and $N$ such that
$$
u(x,t) \geq h_0\quad \hbox{in}\quad {\cal B}^-_{\theta r}.
$$
\end{lemma}
{\it Proof:} We consider $$w=\ln^+(\fr{h}{u+h^{\fr98}}),$$ for
$0<h<1$, to be decided. By applying Lemma 3.3 to $w$ and scaling, we
have
$$ -\!\!\!\!\!\!\int_{{\cal B}^-_{\theta r}}( w-I_0)_+^2 \leq
C(B,\la,N)\fr{\theta \beta r^2}{|{\cal B}^-_{\theta r}|} \int_{{\cal
B}^-_{\beta r} }|Xw|^2.
$$
Let $\tilde{u}={\fr u h}$, then $\tilde{u}$ satisfies the conditions
of Lemma 3.1. We can get similar estimates as (3.2), (3.5), (3.7)
and (3.8), hence we have
$$
\ba{lllllll} && C(B,\la,N)\fr{\theta r^2}{|{\cal B}^-_{\theta r}|}
\int_{{\cal B}^-_{\beta r}}|Xw|^2\\ \\
& &\leq C(B,\la,N)\fr{\theta r^2}{|{\cal B}^-_{\theta r}|}[\fr
{C(B,\la,N)} {\beta^{3Q}(1-\beta)^2} +\fr 45\ln(h^{-\fr 18})]
mes(K_{\beta r}\times S_{\beta r})\\ \\
& &\leq C(\theta,B,N,\la) \ln(h^{-\fr 18}),
 \ea\leqno(3.22)
$$
where $\theta$ has been chosen. By Lemma 2.1, there exists a
constant, still denoted by $\theta$, such that for $z \in {\cal
B}^-_{\theta r}$,
$$
w-I_0\leq C(B,\la,N) (\ln(h^{-\fr 18}))^{\fr 12} .\leqno(3.23)
$$
Therefore we may choose $h_0$ small enough, so that
$$
C (\ln (\fr {1}{h_0^{\fr 18}}))^{\fr 12}\leq  \ln (\fr {1}{2h_0^{\fr
18}})-\lambda_0\ln (\fr {1}{h_0^{\fr 18}}),
$$
then (3.16) and (3.23) gives
$$
\sup_{{\cal B}^-_{\theta r}}\fr{h_0}{u+h_0^{\fr 98}}\leq \fr
{1}{2h_0^{\fr 18}},
$$
which implies $\inf_{{\cal B}^-_{\theta r}}u\geq h_0^{\fr 98}$, then
we have finished our proof.\\

{\bf Proof of Theorem 1.1.} We may assume that $M={\sup}_{{\cal
B}^-_{r}}(+u)={\sup}_{{\cal B}^-_{r}}(-u)$, otherwise we replace $u$
by $u-c$, since $u$ is bounded locally. Then either $1+\fr u M$ or
$1-\fr u M$ satisfies the assumption of Lemma 3.5, and we suppose
$1+\fr u M$ does it, thus Lemma 3.5 implies existing $h_0>0$ such
than $\inf_{{\cal B}^-_{\theta r}}(1+\fr u M)\geq h_0,\, $ i.e.
$u\geq M(h_0-1)$,\,then
$$
Osc_{{\cal B}^-_{\theta r}}u\leq M-M(h_0-1)\leq
(1-\fr{h_0}{2})Osc_{{\cal B}^-_{r}}u,
$$
which implies the $C^{\a}$ regularity of $u$ near point $(0,0)$ by
the standard iteration arguments. By the left invariant
translation group action, we know that $u$ is $C^{\a}$ in the
interior.


\begin{thebibliography}{WWW}
\bibitem [1] {BLU} A. Bonfiglioli, E. Lanconelli and F. Uguzzoni, {\it Stratified Lie Groups and Potential
Theory for their Sub-Laplacians}, Springer-Verlag Berlin Heidelberg,
2007

\bibitem [2]{B1} M. Bramanti, M. C. Cerutti and M. Manfredini.
{\it $L^p$ estimates for some ultraparabolic equations}, J. Math.
Anal. Appl., 200 (2) 332-354 (1996).

\bibitem [3]{C1}C. Cinti, A. Pascucci and S. Polidoro, {\it Pointwise estimates
for solutions to a class of non-homogenous Kolmogorov equations},
Mathematische Annalen, Volume 340, n.2, pp.237-264, 2008

\bibitem [4]{CP}C. Cinti and S. Polidoro, {\it Pointwise local estimates and
Gaussian upper bounds for a class of uniformly subelliptic
ultraparabolic operators}, J. Math. Anal. Appl. 338 (2008) 946-969


\bibitem [5] {G1} G. B. Folland, {\it Subellitic estimates and
function space on nilpotent Lie groups}, Ark. Math., 13 (2):
161-207, (1975).

\bibitem [6] {FP}  M. Di Francesco and S. Polidoro, {\it Harnack inequality for a class of degenerate parabolic equations of
Kolmogorov type}. Adv. Diff. Equ. 11, 1261¨C1320 (2006)

\bibitem [7] {HK}P. Hajlasz and P. Koskela, {\it Sobolev met Poincar$\acute{e}$}, Mem. Amer. Math.
Soc. 145 (2000) x+101.

\bibitem [8]{EL}A. E. Kogoj and E. Lanconelli, {\it An invariant Harnack inequality for a
class of hypoelliptic ultraparabolic equations}, Mediterr. J. Math.
1 (2004) 51¨C80.

\bibitem [9] {S1} S. N. Kruzhkov, {\it A priori bounds and some properties
of solutions of elliptic and parabolic equations}, Math. Sb.
(N.S.) 65 (109) 522-570, (1964).
\bibitem [10] {S2} S. N. Kruzhkov, {\it A priori bounds for generalized
solutions of second-order elliptic and parabolic equations},
(Russian) Dokl. Akad. Nauk SSSR 150  748--751, (1963).
\bibitem [11] {L1} E. Lanconelli and S. Polidoro, {\it On a class of
hypoelliptic evolution operaters}, Rend. Sem. Mat. Univ. Politec.
Torino, 52,1 (1994), 29-63, 1994

\bibitem [12]{L1}A. Lunardi, {\it Schauder estimates for a class
of degenerate elliptic and parabolic operators with unbounded
coefficients in $R^N$}. Ann. Scuola Norm. Sup. Pisa Cl. Sci. 24(4),
133-164 (1997) 23.

\bibitem [13] {M1} M. Manfredini and S. Polidoro, {\it Interior
regularity for weak solutions of ultraparabolic equations in the
divergence form with discontinuous coefficients}, Boll Unione Mat.
Ital. Sez. B Artic. Ric. Mat. (8), 1 (3) 651-675, (1998).

\bibitem [14]{M1}M. Manfredini, {\it The Dirichlet
problem for a class of ultraparabolic equations}. Adv. Diff. Equ. 2,
831-866 (1997) 24.

\bibitem [15] {J1} J. Moser, {\it On Harnack's theorem for elliptic differential
equations}, Comm. Pure Appl. Math. 14  577--591 (1961).

\bibitem [16] {J2} J. Moser, {\it A Harnack inequality for parabolic differential
equations}, Comm. Pure Appl. Math. 17 101--134 (1964).

\bibitem [17] {J0} J. Nash, {\it Continuity of solutions of parabolic and
elliptic equations}, Amer. J. Math., 80, 931-954, (1958).

\bibitem [18] {P1} A. Pascucci and S. Polidoro, {\it The moser's
iterative method for a class of ultraparabolic equations}, Commun.
Contemp. Math. Vol. 6, No. 3 (2004) 395-417.

\bibitem [19] {P2} S. Polidoro and M. A. Ragusa, {\it  H\"older regularity
for solutions of ultraparabolic equations in divergence form},
Potential Anal. 14 no. 4, 341--350  (2001).

\bibitem [19] {WZ}W. Wang and L. Zhang, {\it The $C^{\a}$ regularity of a class of
non-homogeneous ultraparabolic equations}, arXiv:math.AP/0711.3411.

\bibitem [21] {X1} Z. P. Xin and L. Zhang {\it On the global existence of
solutions to the Prandtl's system}, Adv. in Math. 181 88-133
(2004).
\bibitem [22] {X2} Z. P. Xin,  L. Zhang and J. N Zhao, {\it Global well-posedness
for the two dimensional Prandtl's boundary layer equations},
preprint.

\bibitem [23] {Z1} L. Zhang, {\it The $C^\a$ reglarity of a class of
ultraparabolic equations}, arXiv:math.AP/0510405v2 25Dec 2006










\end{thebibliography}
\end{document}